\magnification=\magstep1
\font\T=cmr10 at 17pt
\font\Tt=cmr10 at 14pt
\font\t=cmr10 at 12pt
\def\bigtitle{\T}
\def\mtitle{\Tt}
\def\author{\t}

\footnote{} {This paper was written with the help of the Fonds National
Suisse de la
Recherche Scientifique.} \par

\footnote{} {Mathematics subject classification: 14J17 32S25 57M25}

\centerline{\bigtitle  The boundary of the Milnor fiber}

\vskip.1in

\centerline{\bigtitle of Hirzebruch surface singularities }

\vskip.3in

\centerline{\author  Fran\c coise Michel, Anne Pichon and
Claude Weber }

\vskip.3in

\noindent {\bf Adresses.}

\vskip.1in

\noindent Fran\c coise Michel / Laboratoire de Math\' ematiques Emile Picard  /
Universit\' e Paul Sabatier / 118 route de Narbonne / F-31062 Toulouse
Cedex 04 / FRANCE

e-mail: fmichel@picard.ups-tlse.fr

\vskip.1in

\noindent Anne Pichon / Institut de Math\' ematiques de Luminy / UPR 9016
CNRS / Case
907 / 163 avenue de Luminy / F-13288 Marseille Cedex 9 / FRANCE

e-mail: pichon@iml.univ-mrs.fr

\vskip.1in

\noindent Claude Weber / Section de Math\' ematiques / Universit\' e de
Gen\`eve / CP
64 / CH-1211 Gen\`eve 4 / SUISSE

e-mail: Claude.Weber@math.unige.ch

\vskip.3in

\noindent {\bf Abstract.} We give the first (as far as we know) complete
description of
the boundary of the Milnor fiber for some non-isolated singular germs of
surfaces in
${\bf C}^3$. We study irreducible (i.e.
$gcd~ (m,k,l) = 1$) non-isolated (i.e. $1 \leq k < l$) Hirzebruch hypersurface
singularities in ${\bf C}^3$ given by the equation $z^m - x^ky^l = 0$. We
show that the
boundary $L$ of the Milnor fiber is always a Seifert manifold and we give
an explicit
description of the Seifert structure. From it, we deduce  that:

1) $L$ is never diffeomorphic to the boundary of the normalization.

2) $L$ is a  lens space iff $m = 2$ and $k = 1$.

3) When $L$ is not a lens space, it is never orientation preserving
diffeomorphic to the
boundary of a normal surface singularity.

\vskip.3in

{\mtitle 1. Introduction.}

\vskip.2in

In [MP] the authors prove, among other facts, that the boundary $L$ of the
Milnor fiber
of a non-isolated hypersurface singularity in $\bf{C}^3$ is a Waldhausen
manifold
(non-necessarily "reduziert").

\vskip.1in

In this paper, we apply the general method of [MP] to the study of Hirzebruch
singularities, defined by the equation

$$z^m - x^ky^l = 0$$

We assume that the germ is irreducible, which amounts to ask that
$gcd(m,k,l) = 1$.
We also assume that $1 \leq k \leq l$ to avoid redundancies and that $m
\geq 2$ in order
to have a genuine singularity.

Hirzebruch proved in [H] that the boundary $\tilde L$ of the normalization
is a lens
space and he gave an explicit description of the minimal resolution as a
bamboo-shaped
graph of rational curves. See also [HNK]. We call "bamboo" a connected
graph whose
vertices have at most two neighbours.  We briefly recall this result in
section 2.

\vskip.1in

We prove in theorem 3.1  that $L$ is always a Seifert manifold. Its
canonical star-shaped
plumbing graph is described in theorem 4.2.

\par

When $m \geq 3$ or when $m = 2$ and  $k \geq 2$ the plumbing graph for $L$
is never a
bamboo of rational curves. A little computation shows then (see corollary
4.3) that $L$
is never orientation-preserving diffeomorphic to the boundary of a normal
surface
singularity.

\par

When $m = 2$ and $k = 1$, the plumbing graph is a bamboo of rational
curves. But it is different from the Hirzebruch one. Indeed, the
corresponding lens
spaces do not have the same fundamental group.

\par

In [MP] it is stated that the boundary $L_t$ of the Milnor fiber of a
non-isolated
hypersurface singularity in $\bf{C}^3$ is never diffeomorphic to the boundary
$\tilde{L}_0$ of the normalization. This result is exemplified here in a
very explicit
way, because we are able to compare the two corresponding plumbing graphs
for any
Hirzebruch singularity.

\vskip.1in

The more general case of germs having equation $z^m - g(x,y) = 0$ is
treated in [MPW].
The proofs we present here are self-contained, i.e. independant from [MP]
and from
[MPW].

\vskip.1in

The first named author had the idea to study Hirzebruch singularities while
reading
Egbert Brieskorn beautiful article [B].

\vskip.1in

We thank Walter Neumann for very pleasant discussions during the meeting
and for
attractiong our attention to the computation of the invariant $e_0$ for Seifert
manifolds.

\vskip.3in

{\mtitle 2. Plumbing graphs.}

\vskip.2in

The 3-dimensional manifolds we consider are compact and oriented. In many
cases, they
are oriented as the boundary of a complex surface. To describe these
manifolds, we use
plumbing graphs and we follow [N] as closely as possible. Recall that a
vertex of a
plumbing graph carries two weights: the genus $g$ of the base space and the
Euler
number $e \in \bf{Z}$. In this paper we always have $g \geq 0$~ i.e. the
base surfaces
are orientable. Particuliarly useful are the bamboos for lens spaces and the
star-shaped graphs for "general" Seifert manifolds.

\vskip.1in

The lens space $L(n,q)$ is defined as the quotient of the sphere $S^3
\subset \bf{C}^2$
(oriented as the boundary of the unit 4-ball, equiped with the complex
orientation) by
the action $C_{n,q}$ of the n-th roots of unity given by $\zeta (z_1,z_2) =
(\zeta z_1,
{\zeta}^q z_2)$ with $0 < q < n$ and $gcd(n,q) = 1$. The canonical plumbing
graph for
$L(n,q)$ is the bamboo of rational curves with Euler numbers, from left to
right, $(e_1,
e_2, ... , e_u)$ defined as $e_i = -b_i$. The integers   $b_i$ are defined
by  $b_i \geq
2$ together with

$$\def\un{ 1\hfill} {n \over q} = b_1 -{\un\over b_2 - \displaystyle{\un\over
b_3 - \displaystyle{\un\over \ddots - \displaystyle{\un\over b_u}}}}$$

As in [N], we summarize the continued fraction expansion as $[b_1, b_2, ...
b_u]$.

\vskip.1in

The Seifert manifolds (with unique Seifert foliation) are described by a
star-shaped
graph. See [N] corollary 5.7. All vertices, except possibly the central
one, have genus
zero and Euler number $e \leq -2$.

\vskip.1in

We now consider Hirzebruch singularity $z^m - x^ky^l = 0$. The boundary
$\tilde L$ of
its normalization is the lens space $L(n,q)$ where $n$ and $q$ are computed
as follows.
Let $d_k = gcd (m,k)$ and $d_l = gcd(m,l)$. Then

$$n = {m \over {d_k d_l}}$$

To get $q$ let ${\lambda}_0$ be the smallest integral  positive solution of the
equation

$$\lambda l \equiv -kd_l ~(\bmod m)$$

in the unknown $\lambda$. This solution ${\lambda}_0$ is divisible by $d_k$
and we have

$$q = {{\lambda}_0 \over d_k}$$

\vskip.1in

The special case $d_k = 1 = d_l$ is more pleasant. Then

\vskip.1in

\centerline {$n = m$ and $q = {\lambda}_0$}

\vskip.1in

where ${\lambda}_0$ is the smallest positive solution of the equation
$\lambda l \equiv
-k ~(\bmod m)$. See [BPV].

\vskip.1in

The description we give below in theorem 4.1  for the boundary $L$ of the
Milnor fiber is
in sharp contrast with the classical result (essentially Hirzebruch thesis)
about the
boundary $\tilde L$ of the normalisation. For instance, if $m$ is fixed,
$\tilde L$
depends only on the residue classes $(\bmod m)$ of $k$ and $l$. This is not
the case for
$L$. See section  5   below for an example.

\vskip.3in

{\mtitle 3. Vertical monodromies.}

\vskip.2in

Let $f(x,y,z) = z^m - x^ky^l$ be an irreducible germ (i.e. $gcd(m,k,l) = 1$) of
hypersurface in ${\bf C}^2$ with a singular point (i.e. $2 \leq m$) at the
origin.
Recall that we assume that $1 \leq k \leq l$ to avoid redundancies.

\par

In this paper, we use for technical reasons a polydisc $B(\alpha) =
B^2_{\alpha}
\times B^2_{\alpha} \times B^2_{\epsilon}$ with $0 < \alpha \leq \epsilon$ and
$\alpha^{k+l} < \epsilon^m$ in place of the standard Milnor ball
$B^6_\epsilon = \lbrace
P \in {\bf C}^3 $ with $\vert P \vert \leq \epsilon \rbrace$. The equation
of $f$ being
quasi-homogeneous, for any $B(\alpha)$ there exists  $\eta$ with $0 < \eta
\ll \alpha$
such that the restriction of $f$ on $B(\alpha) \cap f^{-1}(B^2_\eta
\setminus \lbrace 0
\rbrace )$ is a locally trivial fibration on $(B^2_\eta \setminus \lbrace 0
\rbrace)$
and such that this fibration does not depend on $\alpha$ up to isomorphism.
Let $S$ be
the boundary of $B(\alpha)$. The condition $\alpha^{k+l} < \epsilon^m$
implies that we
may choose $\eta$ with  $0 < \eta \ll \alpha$ such that $L_t = f^{-1}(t)
\cap S$
is contained in $\lbrace $ $(x,y,z) \in {\bf}C^3$ such that $\vert x \vert
= \alpha$ or
$\vert y \vert = \alpha$ $\rbrace$ for all $t$ with $0 \leq \vert t \vert
\leq \eta$.
For such a $\eta$, if $t \in B^2_\eta \setminus {0}$ ~we say that $F_t =
B(\alpha) \cap
f^{-1}(t)$ ~is "the" {\bf Milnor fiber} of $f$ and that $L_t = F_t \cap S$
~is "the"
{\bf boundary of the Milnor fiber } of $f$. From now on, we write $L = L_t$
for a
chosen $t$ such that $0 < \vert t \vert \leq \eta$.

\vskip.1in

We will now describe $L$ as the union of $M' = L \cap \lbrace \vert x \vert
= \alpha
\rbrace$ and $M'' = L \cap \lbrace \vert y \vert = \alpha \rbrace$.

\vskip.1in

\noindent {\bf Theorem 3.1}. The boundary $L$ of the Milnor fiber of $z^m -
x^ky^l$ is
a Seifert manifold. Moreover, the projection on the z-axis is constant on
each Seifert
leaf.

\vskip.1in

\noindent {\bf

 Proof of theorem 3.1}. Let $\varphi : M' \rightarrow {\bf C}^3$ be
defined by
$\varphi (x,y,z) = (x,z,f(x,y,z))$. Hence we have $\varphi (M') \subset
S^1_\alpha
\times B^2_\epsilon \times \lbrace t \rbrace$. The singular locus $\Sigma
(f)$ of $f$
satisfies the equation ${\partial f \over \partial y} = 0$ i.e.
$lx^ky^{l-1} = 0$. But
we have $M' \subset {\lbrace \vert x \vert = \alpha} \rbrace$. Hence we
have $\Sigma
(\varphi) = \cup_{i=1}^{m} (S^1_\alpha \times \lbrace 0 \rbrace \times
\lbrace z_i
\rbrace)$ where $z^m_i = t$.

\par

The set of singular values $\Delta (\varphi ) = \varphi ( \Sigma ( \varphi
))$  of the
map  $\varphi$ is the union of the $m$ circles $S^1_\alpha \times \lbrace
z_i \rbrace
\times \lbrace t \rbrace $ where $z^m_i = t$.

\par

We fill $\varphi (M')$ with the circles $S^1_\alpha \times \lbrace c
\rbrace  \times
\lbrace t \rbrace$ where $c \in B^2_\epsilon$ and $\vert c^m - t \vert \leq
\alpha^{k+l}$. As $\Delta ( \varphi )$ is the union of $m$ of these circles
, we
pull-back this (trivial) fibration of $\varphi (M')$
in circles to obtain a Seifert foliation on $M'$. The Seifert leaves are
defined as the
intersection $M' \cap \lbrace z = c\rbrace$.

\par

Replacing $\varphi$ by the restriction to $M''$ of the morphism $ (x,y,z)
\mapsto (y,z,
f(x,y,z))$ we see that, in a symmetric way, the intersections $M'' \cap
\lbrace z = c
\rbrace$ fill $M''$ with a Seifert foliation in circles. The Seifert leaves
of $M'$ and
of $M''$  are defined by the same equation $L \cap \lbrace z = c \rbrace$,
so they
coincide on $T = M' \cap M''$. {\bf End of proof of theorem 3.1.}

\vskip.1in

Let $\pi _x : M' \rightarrow S^1_\alpha$   (resp $\pi _y : M'' \rightarrow
S^1_\alpha$)
be the restriction to $M'$ (resp $M''$) of the projection on the x-axis
(resp the
y-axis). Let $a \in S^1_\alpha$. Now let $G' = {\pi _x}^{-1} (a)$ and $G''
= {\pi
_y}^{-1} (a)$.

\vskip.1in

\noindent {\bf Theorem 3.2}. $\pi _x$ and $\pi _y$ are locally trivial
differentiable
fibrations over $S^1_\alpha$. Moreover:

\noindent 1) The fibers of $\pi _x$ (resp $\pi _y$) are diffeomorphic to
the Milnor
fiber of the plane curve germ $z^m - y^l$ (resp $z^m - x^k$).

\noindent 2) The fibers of $\pi _x$ (resp $\pi _y$) meet transversaly the
Seifert leaves
of $M'$ (resp $M''$) constructed in the proof of theorem 3.1.

\vskip.1in

\noindent {\bf Proof of theorem 3.2}. The singular locus of $\pi _x$ is
defined by
$lx^ky^{l-1} = 0$ and $mz^{m-1} = 0$. But, if $(x,y,z) \in M'$ we have
$\vert x \vert =
\alpha$ and  $z^m - x^ky^l = t$ ~with $0 < \vert t \vert$. So $\pi _x$ has
no singular
point. It is easy to see that the restriction of ${\pi}_x$ to $\partial M'$
is a
submersion onto $S^1_\alpha$. As
$M'$ is a compact differentiable manifold,
$\pi _x$ is a differentiable fibration. The situation is symmetric for $\pi
_y$.

\par

Now, we have chosen $a \in S^1_\alpha$ and $t$ such that  $0 < \vert t
\vert \leq \eta$
where $\eta$ is very small. By definition we have $G' =  \lbrace (a,y,z)
~with~ z^m -
a^ky^l = t ~and~ (y,z) \in S^1_\alpha \times B^2_\epsilon \rbrace$ and also
 $G'' =  \lbrace (x,a,z) ~with~ z^m -  x^ka^l = t ~and~ (x,z) \in
S^1_\alpha \times
B^2_\epsilon \rbrace$. Hence, the assertion 1) is obvious.

\par

To prove 2) let $b$ be any $l^{th}$ root of $(a^{-k}(c^m - t))$ and let $P
= (a,b,c)
\in G'$. The Seifert leaf containing $P$ is parametized by $(e^{i\theta}a,
e^{-i\theta
{k \over l}} b, c)$ with, say, $\theta \in {\bf R}$. Hence, the Seifert
leaves are
oriented and transverse to the hyperplane $H_a = \lbrace x = a \rbrace$ for
all $a \in
S^1_\alpha $. The situation is symmetric for $M''$. {\bf End of proof of
theorem 3.2}.

\vskip.1in

\noindent {\bf Remarks.} 1. If $k = l = 1$ the germ $f$ has an isolated
singular point
at the origin. In this case, theorem 3.2 shows that $G'$ and $G''$ are
discs and that
 $M'$ and $M''$ are solid torii. Hence $L$ is a lens space, diffeomorphic
to $L_0 =
\tilde L$.

\par

2. If we assume that dim$\Sigma (f) = 1$ then we have $l \geq 2$ and the
x-axis $D' =
\lbrace (x,0,0) ~with ~x \in {\bf C} \rbrace $ is a component of $\Sigma
(f)$. Then,
theorem 3.2 implies that $G'$ is never diffeomorphic to a disc and that
$M'$ is not a
solid torus. When $D' \subset \Sigma (f)$ we say in [MP] that $M'$ is the
{\bf vanishing
zone} around $D'$. When $k \geq 2$ then $D'' = \lbrace (0,y,0) ~with ~y \in
{ \bf C}
\rbrace$ is the second component of $\Sigma(f)$ and $M''$ is the vanishing
zone around
$D''$.

\vskip.1in

We now proceed to the definition of the vertical monodromy. Let $h' : G'
\rightarrow
G'$ be the diffeomorphism defined by the first return along the (oriented)
leaves of
$M'$. Theorem 3.2 implies that $h'$ is a monodromy for the fibration $\pi _x$.

\vskip.1in

\noindent {\bf Definition.} We call $h'$ the {\bf vertical monodromy} for
$D'$.

\vskip.1in

Likewise, the first return along the ( oriented) Seifert leaves of $M''$ is a
diffeomorphism
$h'' : G'' \rightarrow G''$. We call it the vertical monodromy for $D''$.

\vskip.1in

In conclusion, we  know that $M'$ is the mapping torus of $h'$ acting on
$G'$ and
that $M''$ is the mapping torus of $h''$ acting on $G''$. We wish now to
describe in
details the vertical monodromies.

\noindent {\bf Notations.} Let $ d = gcd (k,l) ~;~ {\bar l} = {l \over d}
~; ~{\bar k}
= {k \over d} ~; ~ d_l = gcd (m,l) ~; ~ d_k = gcd (m,k) $.

\vskip.1in

\noindent {\bf Remark.} As $f$ is assumed to be irreducible, we have $gcd
(m,k,l) = 1$
and $\bar k$ is prime to $d_l$ (resp $\bar l$ is prime to $d_k$). Moreover,
$G'$ has
$d_l$ boundary components and $G''$ has $d_k$ boundary components.

\vskip.1in

\noindent {\bf Theorem 3.3.}  The vertical monodromy $h'$ (resp $h''$) has
finite
order $\bar l$ (resp $\bar k$). Moreover:

\par

\noindent 1. If $\bar l \geq 2$ (resp $\bar k \geq 2$) then  $h'$ (resp
$h''$) has
exactly
$m$ fixed points and any non-fixed point has order $\bar l$ (resp $\bar k$).

\par

\noindent 2. At each fixed point $h'$ (resp $h''$)  acts locally as a
rotation of angle
$-({{\bar k} \big/ {\bar l}}) 2\pi$ (resp $-({{\bar l} \big/ {\bar k}}) 2\pi$).

\vskip.1in

\noindent {\bf Proof of theorem 3.3.} As in the proof of theorem 3.2, we
consider $P =
(a,b,c) \in G'$. We have seen that $(e^{i\theta}a, e^{-i\theta {k \over l}}
b, c)$ for,
say,
$\theta \in {\bf R}$ is a parametrization of the Seifert leaf which
contains $P$. Hence

$$ (\star ) ~~~~~ h'(P) = (a, e^{-2i\pi {k \over l}}b,c)$$

As ${k \big/ l} = {{\bar k} \big/ {\bar l}}$ with $\bar k$ prime to $\bar
l$, we see
that $h'$ has order $\bar l$ on each $P = (a,b,c)$ with $b \neq 0$.

Then, if $\bar l \geq 2$,  it is clear that $h'(P) = P$ iff $ b = 0 $. Then
$c^m = t$
and $h'$ has exactly
$m$ fixed points, i.e. the points $ \lbrace (a,0,z_i) \rbrace$  where
$z^m_i = t$.

The formula $( \star )$ implies directly the last statement of theorem 3.3.
{\bf End of proof of theorem 3.3}.

\vskip.1in

\noindent {\bf Corollary 3.4.} The intersection $T = M' \cap M''$ is a torus.

\vskip.1in

\noindent {\bf Proof of corollary 3.4.} Indeed, $T$ is the mapping torus of
$h'$ acting
on the $d_l$ boundary components of $G'$. As $\bar k$ is prime to $\bar l$
~the formula
$(
\star )$ in the proof of theorem 3.3 implies that $h'$ permutes
transitively the
boundary components of $G'$. {\bf End of proof of corollary 3.4.}

\vskip.1in

\noindent {\bf Remark.} $G'$, $G''$ and $ F_t =  f^{-1} (t)  \cap B$ are
oriented by
the complex structure. $L$ is oriented as the boundary of $F_t$ and this
orientation
induces one on $M'$ and $M''$.

\vskip.1in

\noindent {\bf Theorem 3.5.} Orient $T = M' \cap M''$ as the boundary of
$M''$. Orient
$\partial G'$ (resp $\partial G''$) as the boundary of $G'$ (resp $G''$).
Then the
intersection number on $T$ of $\partial G'$ with $\partial G''$ is equal to
$-m$.

\vskip.1in

\noindent {\bf Proof of theorem 3.5.} Let $\pi : L \rightarrow B^2_\alpha
\times
B^2_\alpha$
 be the restriction on $L$ of the projection $(x,y,z) \mapsto (x,y)$. The
restriction
of $\pi $ to $ T = M' \cap M''$ is a regular covering of order $m$.
Moreover, we have
$\pi (G') = \lbrace a \rbrace \times S^1_\alpha$ and $\pi (G'') =
S^1_\alpha \times
\lbrace a \rbrace$. The complex structure of ${\bf C}^2$ induces an
orientation on
$B^2_\alpha \times B^2_\alpha$. Let $S^1_\alpha \times S^1_\alpha = \pi (T)$ be
oriented as the boundary of $B^2_\alpha \times S^1_\alpha$. The
intersection number of
$\lbrace a \rbrace \times S^1_\alpha$ with $S^1_\alpha \times \lbrace a
\rbrace$ in
$S^1_\alpha \times S^1_\alpha$ is equal to $(-1)$.
The covering projection $\pi$ being compatible with orientations, this
proves that the
intersection number we are looking for is equal to $(-m)$.
{\bf End of proof of theorem 3.5.}

\vskip.3in

{\mtitle 4. The Seifert structure on the boundary of the Milnor fiber.}

\vskip.2in

\noindent  {\bf Theorem 4.1.} The Seifert invariants (associated to the
Seifert structure
described in section 3) for the boundary
$L$ of the Milnor fiber of a Hirzebruch singularity are as follows:

1. The genus $g$ of the base space is equal to $(m-1)(d-1)$ where $d = gcd
(k,l)$.

2. The integral Euler number $e$  is equal to $m$.

3. Let $ {\bar l} = {l \over d}$ and $ {\bar k} = {k \over d}$. Then $L$
has $2m$
(possibly) exceptional leaves.

There are $m$ of them with Seifert invariants $(\alpha ', \beta ')$ defined
by $\alpha
' = \bar l$ and $\beta '$ given by $(-\bar k) \beta ' \equiv 1 \bmod {\bar
l}$ and $0 <
\beta ' < \bar l$ in normalized form.

There are $m$ of them with Seifert invariants $(\alpha '', \beta '')$
defined by
$\alpha '' = \bar k$ and $\beta ''$ given  by $(-\bar l) \beta '' \equiv 1
\bmod\bar
k$ and $0 < \beta '' < \bar k$.

\vskip.1in

\noindent {\bf Comments.} 1. The singularity is isolated iff $k = l = 1$.
Of course in
this case we have ${\tilde L} = L$. The theorem above says that $L$ has no
exceptional
leaf, that $g = 0$ and that $e = m$. Hence $L$ is the lens space
$L(m,m-1)$. We are
happy to see that this agrees with Hirzebruch result.

\noindent Assume from now on that $1 \leq k$ and that $2 \leq l$.

\noindent 2. Under this hypothesis $L$ is a lens space iff $m = 2$ and $k =
1$. (Quick
proof: To get a lens space we need $g = 0$ and the theorem says that this
is equivalent
to $d = 1$. Then we can admit at most two exceptional leaves. Hence $k = 1
$ and $m =
2$). The lens space is $L(2l, 1)$. On the other hand ${\tilde L} = L(1,1) =
S^3$ when
$l$ is even and ${\tilde L} = L(2,1) = P^3 (\bf R)$ when $l$ is odd.

\noindent 3. If $3 \leq m$ or if $m = 2$ and $2 \leq k$ then at least one
of the two
following statements is true:

i) $g$ is strictly positive

ii) $L$ has strictly more than two exceptional leaves.

\vskip.1in

We describe the canonical plumbing graph in the next theorem. Its proof follows
immediately from theorem 4.1 and from the recipes in [N].

\vskip.1in

\noindent {\bf Theorem 4.2.} 1. If $k = l = 1$ the canonical plumbing graph
is a bamboo
of rational curves, having $(m-1)$ vertices with Euler number equal to
$(-2)$. This is
the singularity $A_{m-1}$.

\noindent Assume from now on that $1 \leq k$ and that $2 \leq l$.

\noindent 2. If $k = 1$ and $m = 2$ the plumbing graph has just one vertex
with $g = 0$
and $e = -2l$.

\noindent 3. Assume either that $3 \leq m$ or that $m = 2$ and $2 \leq k$.
Then the
canonical plumbing graph is never a bamboo of rational curves. More precisely:

\noindent 3a. If $k = l$ the graph has just one vertex with $g =
(m-1)(d-1)$ and $e =
m$. Notice that $g$ is strictly positive because $d = k = l > 1$.

\noindent 3b. If $k$ divides $l$ but $ k \neq l$ the graph is star-shaped
with $m$
branches. The central vertex has $g = (m-1)(d-1)$ and $e = 0$. Each branch
has just one
vertex (tied to the central vertex by an edge). Its weights are $ g = 0$
and $e = -{ l
\over k}$.

\noindent 3c. If $k$ does not divide $l$ then the graph is star-shaped with
$2m$
branches. The central vertex has $g = (m-1)(d-1)$ and $ e = -m$.

There are $m$ branches which are a bamboo of rational curves with $e'_i =
-b'_i$ and
$b'_i$ defined by $b'_i \geq 2$ and

$${\alpha ' \over {\alpha ' - \beta '}} = [b'_1, ... , b'_u]$$

The vertex carrying the number $1$ is joined to the central vertex by an edge.

There are also  $m$ branches which are a bamboo of rational curves with
$e''_i = -b''_i$
and
$b''_i$ defined by $b''_i \geq 2$ and

$${\alpha '' \over {\alpha '' - \beta ''}} = [b''_1, ... , b''_v]$$

Again, the vertex carrying the number $1$ is joined to the central vertex
by an edge.

\vskip.1in

\noindent {\bf Corollary 4.3.} If $L$ is not a lens space, it is never
orientation
preserving diffeomorphic to the boundary of a normal surface singularity.

\vskip.1in

\noindent {\bf Proof of corollary 4.3.} $L$ is not a lens space iff we are
in case 3.
We claim that the intersection form associated to the canonical plumbing
graph is never
negative definite. In cases 3a and 3b this is obvious since the
self-intersection of the
central vertex is $\geq 0$.

Let us suppose that we are in case 3c. We compute the rational Euler number
$e_0$ of
the Seifert structure on $L$. By definition

$$e_0 = e - \sum {{\beta _i} \over {\alpha _i}}$$

>From theorem 4.1 we deduce that

$$e_0 = m - m{\beta ' \over \bar l} -m{\beta '' \over \bar k} $$

Hence:

$${\bar k} {\bar l} e_0 = m({\bar k} {\bar l} -\beta ' \bar k -\beta ''
\bar l)$$

We shall prove later in this section that $({\bar k} {\bar l} -\beta ' \bar
k -\beta ''
\bar l) = 1$. See  lemma 4.6.

Hence

$$e_0 = {m \over {\bar k}{\bar l}} > 0$$

The conclusion follows from [N] Corallary 6 p.300. {\bf End of proof of
corollary 4.3.}

\vskip.1in

\noindent {\bf Proof of theorem 4.1.} We shall compute the Seifert
invariants from the
data provided by the theorems proved in section 3.

\vskip.1in

We first determine the genus $g$. The Euler characteristic $\chi (G')$ is
equal to $(-ml
+m +l)$. The classical formula for ramified coverings implies that the Euler
characteristic $\chi '$ of the quotient of $G'$ by the action generated by
$h'$ is
equal to $(-md + d +m)$. An analogous computation shows that $\chi '' =
\chi '$. Hence
the Euler characteristic $\chi$ of the base space of the Seifert foliation
is equal to
$2(-md + d +m)$ and we get $ g = (m-1)(d-1)$.

\vskip.1in

The computation of the Seifert invariants $(\alpha, \beta)$ is routine if
we use the
dictionary which translates Nielsen invariants into Seifert's.

It is sufficient for us
to consider the following special case. Suppose that the angle of rotation
at a fixed
point of a monodromy $h$ of finite order acting on an oriented surface is
equal to
${{\omega} \over {\lambda}} 2\pi$ with $gcd (\omega , \lambda) = 1$. Define
$\sigma$ as
the integer which satisfies $0 < \sigma < \lambda$ and $\omega \sigma
\equiv 1 ~  (\bmod
\lambda)$. In the mapping torus of $h$, the Seifert invariant $(\alpha ,
\beta )$ for
the exceptional leaf which corresponds to the fixed point is given by $\alpha =
\lambda$ and $\beta = \sigma$ in normalized form. See [M]. The result
follows now
immediately from theorem 3.3.

\vskip.1in

The delicate part of the proof is to determine the Euler number $e$. As we
feel that
this invariant is rather elusive, we prefer to deal with closed objects.

Let $\hat G'$ be the closed surface obtained from $G'$ by attaching a disc
on each of
its $d_l = gcd (m,l)$ boundary components. We have seen (in the proof of
Corollary 3.4)
that the monodromy
$h'$ permutes them transitively. Let $\hat h'$ be "the" finite order
extension  of $h'$
on
$\hat G'$. There is exactly one orbit of $\hat h'$ which corresponds to the
center of
these discs. Its Nielsen invariant $\sigma \big/ \bar l$ is given by

$${\sigma \over \bar l} \equiv -m {\beta ' \over \bar l} ~~  in~~ {\bf Q}
\bmod {\bf
Z}$$

because the sum of all Nielsen quotients is equal to zero in ${\bf Q} \bmod
{\bf Z}$
for a closed surface.

Let $\hat M'$ be the mapping torus of $\hat h'$ acting on $\hat G'$. It is
a closed
Seifert manifold. It has $m$ exceptional leaves with Seifert invariant
$(\alpha ',
\beta ')$ and one with Seifert invariant $(\hat \alpha ', \hat \beta ')$
which we choose
to be defined as

$${\hat \beta ' \over \hat \alpha ' } = -m {\beta ' \over \alpha ' }$$

where $\hat \beta '$ and $\hat \alpha '$ are by necessity chosen to be
relatively
prime. This choice has the advantage that the Euler number $\hat e'$ for
$\hat M' $ is
equal to zero, because the rational Euler number for $\hat M'$ is equal to
zero, as
$\hat M'$ is  the mapping torus of a finite order monodromy acting on a
closed surface.
See  [P].

We proceed along the same path with $G''$ and $h''$ to get a closed Seifert
manifold
$\hat M''$ with analogously defined Seifert invariants.

\vskip.1in

We now state a lemma about glueings of Seifert manifolds. The statement is
painful
(sorry!).

\vskip.1in

\noindent {\bf Lemma 4.4.} Let $V'$ and $V''$ be two closed oriented
Seifert manifolds.
Let
$H'_0$ be a leaf in $V'$ and let $H''_0$ be one in $V''$. Let $N'$ be a
foliated closed
tubular neighborhood of $H'_0$ in $V'$ and let $N''$ be one for $H''_0$ in
$V''$.

Let $s'$ be a section in $V'$
(as usual possibly outside some discs in the base space) giving rise to an
Euler number
$e'$ for $V'$ and a Seifert invariant $(a',b')$ for $H'_0$. In a similar
manner, let
$s''$ be a section in $V''$ giving rise to the Euler number $e''$ for $V''$
and to the
Seifert invariant $(a'',b'')$ for $H''_0$.

Let $\check V' = V' \setminus Int(N')$ and $\check V'' = V'' \setminus
Int(N'')$. Let
$V$ be such that $V = \check V' \cup \check V''$ and $\check V' \cap \check
V'' =
\partial {\check V'} \cap \partial \check V''$. This intersection is a
torus  and we
write $T$ for it. Suppose that the leaves $H'$ from $V'$ and $H''$ from
$V''$ coincide
on $T$ (hence $V$ is Seifert foliated).

Let $m'$ be a meridian for $N'$ on $T$ and let $m''$ be one for $N''$. Let
$IN(m',m'')$
be the intersection number of $m'$ and $m''$ on $T$, where $T$ is oriented
as the
boundary of $\check V''$.

Then the Euler number $e$ for $V$ (corresponding to a section $s$
essentially built from
$s'$ and $s''$) is given by the equality $e = e' + e'' + \bar e$ where
$\bar e$ is
computed from the equation

$$IN(m',m'') = a'b'' + a''b' + a'a'' \bar e$$.

\noindent {\bf Proof of lemma 4.4.} As the section $s$ is built from $s'$
and $s''$ it
follows from the definition of the Euler number  as an obstruction
(evaluated on a
fundamental cycle) that $e$ is the sum of $e'$ and $e''$ plus a
contribution coming
from the fact that $s'$ and $s''$ do not necessarily match along the torus
$T$. The
formula of theorem 3.5 will determine that contribution.

Following Seifert conventions we have

$$m' = a's' + b'H' ~~with~~ a' > 0 ~~~and~~~ m'' = a''s'' + b''H'' ~~with~~
a'' > 0$$

By hypothesis, we have $H' = H'' = H$. Let us choose an orientation
(arbitrarily) for
$H$. From Seifert conventions, this choice orients $s'$ and $s''$ via
$IN(s',H) = +1$
on $T$ oriented as $\partial N'$ and $IN(s'',H) = +1$ on $T$ oriented as
$\partial
N''$. This orients $m'$ on $T = \partial N'$ via $a' > 0$ and $m''$ on $T =
\partial
N''$ via $a'' > 0$.

Notice that a change of orientation of $H$ induces a change of orientation
on both $m'$
and $m''$ and hence the intersection number $IN(m',m'')$ does not change.
Let us
compute that intersection number.

\vskip.1in

$IN(m',m'') = IN((a's' + b'H),(a''s'' + b''H))$

 $= a'a''IN(s',s'') + a'b''IN(s', H) + a''b'IN(H,s'') + b'b''IN(H,H)$

\vskip.1in

We have:

\vskip.1in

\noindent 1) $IN(H,H) = 0$ because the intersection form is alternating.

\noindent 2) $IN(s',H) = + 1$ from Seifert conventions, because $T$ is
oriented as the
boundary of $\check V''$ which is the same as being oriented as the
boundary of $N'$.

\noindent 3) $IN(H,s'') = +1 $  because $IN(s'',H) = +1$  if $T$
is oriented as the boundary of $N''$ and two sign changes occur from the
last equality
to get the first one.

\noindent 4) $IN(s',s'') = \bar e$. To see that the sign is correct, one
way to argue
is to go back to the definition of Euler numbers. Another way is to remark
that this is
the good sign in order to be sure that the sum $e' + e'' + \bar e$ remains
constant
under changes of $s'$ (or $s'')$ near the fiber $H'_0$ (or $H''_0)$.

\vskip.1in

\centerline {{\bf End of proof of lemma 4.4.}}

\vskip.1in

We now use lemma 4.4 to complete the determination of $e$. To make the
argument simpler
let us assume that

$(d_k = gcd (m,k) = 1 ~~; ~~ d_l = gcd (m,l) = 1 ~~; ~~d = gcd (k,l) = 1)$

\vskip.1in

Recall that in this case $\hat M'$ has $m$ exceptional leaves with Seifert
invariant
$\alpha ' = l$ and $\beta '$ defined by $0 < \beta ' < l$ and $(-k) \beta '
\equiv 1
~~(\bmod l)$. $\hat M'$ has one more exceptional leaf with Seifert
invariant $(\hat
\alpha ', \hat \beta ')$ defined by

$${{\hat \beta '} \over {\hat \alpha '}} = -m {{\beta '} \over {l}}$$

As $gcd (m,l) = 1$ we have that $\hat \alpha ' = l$. We have already seen
that $e' = 0$.

\vskip.1in

Similarly, $\hat M''$ has $m$ exceptional leaves with invariant $\alpha ''
= k$ and
$\beta ''$ defined by $0 < \beta '' < k$ and $(-l) \beta '' \equiv 1
~~(\bmod k)$.
$\hat M''$ has one more exceptional leaf with invariant $(\hat \alpha '',
\hat \beta
'')$ defined by

$${{\hat \beta ''} \over {\hat \alpha ''}} = -m {{\beta ''} \over {k}}$$

We have $\hat \alpha '' = k$ because $gcd (m,k) = 1$ and $e'' = 0$.

\vskip.1in

As $gcd (m,l) = 1$ the boundary $\partial G'$ is connected and $\partial
G''$ is
connected because $gcd (m,k) = 1$. As a consequence, the intersection
number $IN
(\partial G', \partial G'') $ is equal to $IN(m',m'')$~~ UP TO SIGN.

\vskip.1in

\noindent {Lemma 4.5.} We have the equality  $IN(m',m'') = - IN( \partial
G', \partial
G'')$.

\vskip.1in

\noindent {\bf Proof of lemma 4.5.} The result comes from a comparison
between the
orientation of meridians coming from Seifert conventions and the
orientation coming
from $\partial G'$ (or $\partial G'')$. What happens is that for one
meridian both
orientations agree and that for the other one they disagree. Which one it
is depends on
the orientation selected for $H$. {\bf End of proof of lemma 4.5.}

\vskip.1in

We go on with the determination of the Euler number. The formula

$$IN(m',m'') = a'b'' + a''b' + a'a''\bar e$$ of lemma 4.4 translates into

$$m = l(-m \beta '') + k(-m \beta ') + kl \bar e$$

Hence we have

$$ (\dagger) ~~~   m(1 + l\beta '' + k\beta ') = kl \bar e$$

\vskip.1in

\noindent {\bf Lemma 4.6.} We have the equality: $(\star \star)$ $1 + l
\beta '' + k
\beta ' = kl$.

\vskip.1in

>From lemma 4.6 and formula ($\dagger$) we deduce that $\bar e = m$ and
hence that $e =
m$~ because $e' = 0 = e''$. This completes the computation of $e$.

\vskip.1in

\noindent {\bf Proof of lemma 4.6.} By definition we have

$$l \beta '' \equiv -1 ~~(\bmod k) ~~and~~ k \beta ' \equiv -1 ~~(\bmod l)$$

Because $gcd(k,l) = 1$ we deduce that

$$l \beta '' + k \beta ' \equiv -1 ~~(\bmod kl)$$

In other words there exists an integer $q$ such that

$$1 + l \beta '' + k \beta ' = qkl$$

As $0 < \beta ' < l$ and $0 < \beta '' < k$ the only possibility is $q =
1$. {\bf End
of proof of lemma 4.6.}

\vskip.1in

By carefully dividing by adequate gcd's
an analogous argument works without assuming that
($d_k = gcd (m,k) = 1 ~~; ~~ d_l = gcd (m,l) = 1 ~~; ~~d = gcd (k,l) = 1)$.
{\bf End of proof of theorem 4.1.}

\vskip1in

{\mtitle 5. Examples.}

\vskip.2in

\noindent {\bf Example 1.} Let us consider the Hirzebruch singularity $z^{12} -
x^5y^{11} = 0$

\vskip.1in

The boundary $\tilde L$ of the normalization is the lens space $L(12,5)$.
Its plumbing
graph is a bamboo of three rational curves with Euler numbers successively
$\lbrace
-3,-2,-3
\rbrace$.

\vskip.1in

The Seifert structure of the boundary $L$ of the Milnor fiber is as follows:

\noindent ($ g = 0$ and $e = 12$). $L$ has 24 exceptional leaves. There are
12 of them
with Seifert invariant
$(\alpha = 11,~ \beta = 2)$ and 12 of them with Seifert invariant $(\alpha
= 5,~ \beta
= 4)$.

\vskip.1in

The plumbing graph of $L$ is star-shaped. The central vertex has weights $g
= 0$ and $e
= -12$. There are 24 bamboos of rational curves  attached to the central
vertex.
Among them,  12  have Euler numbers equal successively to $\lbrace -2, -2,
-2, -2,
-3
\rbrace$ and 12 of them have just one vertex with Euler number equal  to
$\lbrace -5
\rbrace$.

\vskip.1in

\noindent {\bf Example 2.} Let us consider the Hirzebruch singularity $z^{12} -
x^{17}y^{11} = 0$. In order to make the comparison between examples 1 and 2
easier, we
drop the restriction $k \leq l$.

\vskip.1in

The boundary $\tilde L$ of the normalization is the same as in example 1,
because
5 is  congruent to 17 (mod 12).

\vskip.1in

But the boundaries $L$ of the Milnor fibers
are different. In fact, the Seifert invariants for the exceptional leaves
differ. $L$
has 12 leaves with Seifert invariant $(\alpha = 11,~ \beta = 9)$ and 12
leaves with
Seifert invariant $(\alpha = 17,~ \beta = 3)$.

\vskip.1in

The plumbing graph of $L$ is again star-shaped, as it should be. The
central vertex has
again weights $g = 0$ and $e = -12$. There are 24 bamboos of rational
curves  attached
to the  central vertex. Among them, 12
have Euler numbers equal successively to $\lbrace -6, -2 \rbrace$ and 12 of
them have Euler numbers successively equal   to $\lbrace -2, -2, -2, -2, -3, -2
\rbrace$.

\vskip.3in

\noindent {\mtitle 6. Bibliography.}

\vskip.2in

[B] E. Brieskorn: "Singularities in the work of Friedrich Hirzebruch".
Surv. Diff.
Geom. VII, Int. Press, Sommerville, MA (2000), 17-60.

[BPV] W. Barth, C. Peters, A. Van de Ven: "Compact Complex Surfaces".
Ergebnisse
der Mathematik und ihrer Grenzgebiete {\bf 3}, Band 4, Springer Verlag (1984).

[H] F. Hirzebruch: "\"Uber vierdimensionale Riemannsche Fl\"achen
mehrdeutiger analytischeer Funktionen von zwei Ver\"anderlichen". Math.
Ann. {\bf
126} (1953), 1-22.

[HNK] F. Hirzebruch, W. D. Neumann, S. S. Koh: "Differentiable manifolds
and quadratic
forms". Math. Lecture Notes, vol 4, Dekker, New-York (1972).

[J] H. Jung: "Darstellung der Funktionen eines algebraischen K\"orpers zweier
unabh\"angigen Ver\"anderlichen $(x,y)$ in der Umgebung einer Stelle
$(x-a,y-b)$". Jour. reine u. angew. Mathematik {\bf 133} (1908), 289-314.

[M] J. Montesinos: "Classical tessellations and three-manifolds".
Universitext, Sprin-ger Verlag, Berlin (1987).

[MP] F. Michel, A. Pichon: "On the boundary of the Milnor fiber of non-isolated
singularities". IMRN {\bf 43} (2003), 2305-2311.

[MPW] F. Michel, A. Pichon, C. Weber: "An explicit description of the
boundary of the
Milnor fiber for some non-isolated surface singularities in ${ \bf C}^3$".
Manuscript in
preparation.

[N] W. D. Neumann: "A calculus for plumbing applied to the topology of
complex surface
singularities and degenerating complex curves". Trans. AMS {\bf 268}
(1981), 299-344.

[P] A. Pichon: "Fibrations sur le cercle et surfaces complexes". Ann.
Inst. Fourier (Grenoble) {\bf 51} (2001), 337-374.

\vskip.3in

Gen\`eve, le 30 juin 2005.

\end